\documentclass[11pt,  twoside]{amsart}

\usepackage{color}

\usepackage{graphicx}

\usepackage[
hypertexnames=false, colorlinks, citecolor=red, linkcolor=red]{hyperref} 
\hypersetup{bookmarksdepth=3}

\usepackage{geometry}
\usepackage{amsmath,amssymb}
\usepackage{
mathrsfs}
\usepackage{esint}

\allowdisplaybreaks[1]

\newtheorem{thm}{Theorem}[section]
\newtheorem{lm}[thm]{Lemma}

\newtheorem{prop}[thm]{Proposition}

\theoremstyle{definition}

\newtheorem*{df*}{Definition}

\theoremstyle{remark}

\newtheorem*{rem*}{Remark}

\numberwithin{equation}{section}

\newcommand{\ci}[1]{_{ {}_{\scriptstyle #1}}}

\newcommand{\1}{\mathbf{1}}

\newcommand{\la}{\lambda}

\newcommand{\cG}{\mathcal{G}}
\newcommand{\cX}{\mathcal{X}}

\newcommand{\cz}{Calder\'{o}n--Zygmund\ }

\newcommand{\Z}{\mathbb{Z}}

\newcommand{\R}{\mathbb{R}}

\newcommand{\wt}{\widetilde}

\newcommand{\cD}{\mathscr{D}}

\newcommand{\cE}{\mathcal{E}}
\newcommand{\cL}{\mathcal{L}}

\newcommand{\bal}{{\boldsymbol\alpha}}

\newcommand{\fA}{\mathfrak{A}}

\def\cyr{\fontencoding{OT2}\fontfamily{wncyr}\selectfont}
\DeclareTextFontCommand{\textcyr}{\cyr}

{\end{list}}      

\renewcommand{\labelenumi}{(\roman{enumi})}

\newcounter{vremennyj}

\newcommand\cond[1]{\setcounter{vremennyj}{\theenumi}\setcounter{enumi}{#1}\labelenumi\setcounter{enumi}{\thevremennyj}}

\let\oldtocsection=\tocsection

\let\oldtocsubsection=\tocsubsection

\let\oldtocsubsubsection=\tocsubsubsection

\renewcommand{\tocsection}[2]{\hspace{0em}\oldtocsection{#1}{#2}}
\renewcommand{\tocsubsection}[2]{\hspace{1em}\oldtocsubsection{#1}{#2}}
\renewcommand{\tocsubsubsection}[2]{\hspace{2em}\oldtocsubsubsection{#1}{#2}}

\begin{document}

\title
{A remark on two weight 
estimates for positive dyadic operators. }

\author{Sergei Treil}
\thanks{Supported  by the National Science Foundation under the grant  DMS-0800876. }
\address{Dept. of Mathematics, Brown University,   
151 Thayer
Str./Box 1917,      
 Providence, RI  02912, USA }
\email{treil@math.brown.edu}
\urladdr{http://www.math.brown.edu/\~{}treil}

\makeatletter
\@namedef{subjclassname@2010}{
  \textup{2010} Mathematics Subject Classification}
\makeatother




%
%

\date{}

\begin{abstract}
We give a simple proof of the Sawyer type characterization of the two weigh estimate for positive dyadic operators  (also known as the bilinear embedding theorem)
\end{abstract}

\maketitle

\setcounter{tocdepth}{1}
\tableofcontents

%

\section{Introduction}
\label{intro}

The paper deals with the two weight estimates for the so-called positive dyadic operators $T_\alpha$, $\alpha=\{\alpha\ci I\}\ci{I\in\cD}$, $\alpha\ci I \in[0,\infty)$, 
\[
T_\alpha f := \sum_{I\in\cD} \alpha\ci I \left(\int_I f \,d\mu \right) \1\ci I, 
\] 
where $\cD$ is a dyadic lattice in $\R^d$. 

We are presenting a simple proof of the (well-known) fact that the so-called Sawyer type conditions are sufficient for the boundedness of the operator $T_\alpha : L^p(\mu) \to L^p (\nu)$.   The Sawyer type condition essentially mean that the operator $T_\alpha$ and its formal adjoint are uniformly bounded on functions $\1\ci I$, $I\in\cD$ (they are in fact a bit weaker, see exact statement in Theorem \ref{t-main} below), so they are trivially necessary. 

The conditions are named after E.~Sawyer, who in \cite{Sawyer_2weight-max_1982} who proved that such conditions for the maximal operator $M$ are sufficient for its boundedness. Note that here is sufficient to use only the conditions for $M$: no condition for the adjoint operator (which is not well defined, by the way) is needed. 

Later he proved in \cite{Sawyer_2weightFract-1988} that such conditions (now for the operator and its formal adjoint) are sufficient for the two weigh estimates for a wide class of integral operators with non-negative kernels, in particular for fractional integrals and Poisson integrals. 

Note, that while the conditions in \cite{Sawyer_2weight-max_1982} and \cite{Sawyer_2weightFract-1988} look differently from the conditions in Theorem \ref{t-main} below, they can be transformed to the form presented here by the standard ``change of measure'' argument. 

While it was expected that Theorem \ref{t-main} should be true, formally it did not follow from the result in \cite{Sawyer_2weightFract-1988}.  

Theorem \ref{t-main} was first proved for $p=2$ in \cite{NTV-2w} by the Bellman function method.%
\footnote{The author is also familiar with a manuscript by F.~Nazarov,  dated back to the same time as  \cite{NTV-2w}, where this result was proved for all$p\in (1, \infty)$, again using the Bellman function method. However,  this manuscript was never published.} It was proved there under the name of ``bilinear embedding theorem'' and it was stated in the same form as the result in \cite{Sawyer_2weightFract-1988}.

Later in \cite{Lacey-2weiight_discr_pos-2009} this theorem was proved in full generality; in fact the case of operators $L^p(\mu)\to L^q(\nu)$ was also treated there. The construction was quite complicated  but was based on the standard techniques of the modern harmonic analysis like stopping moment, corona decomposition etc. 

Here we present a simple proof of Theorem \ref{t-main}. To make the presentation more clear we are considering here only the hardest case $p=q$. 

There are  two main ideas in the proof. The first one that goes back to F.~Nazarov and was already used in \cite{NTV-2w} is that the testing condition for the operator and its adjoint each are responsible for the bounds on part of the sum in the bilinear embedding: the splitting of the sum is determined by the condition \eqref{split-cond1} below. So the main theorem is reduced to Proposition \ref{p-BE_half}. 

In the proof of Proposition \ref{p-BE_half} the sum is represented as sum of the integrals over stopping cubes, and each integral is split into two parts. The first parts have disjoint supports and the estimate follows from the Carleson Embedding Theorem. 

The second (hard) parts do not have disjoint supports, but the supports are small and satisfy the Carleson measure condition. And the estimate for these parts is obtained by noticing that at some places one can repace the function $g$ by its averages and then estimating  the  averages of $g$ by the averages of $f$ using the splitting condition \eqref{split-cond1} and thus again reducing the bilinear estimate to the Carleson embedding theorem for $f$. This replacement of averages of $g$ by the averages of $f$ is the second idea of the proof.

\section{
The bilinear embedding theorem}

\begin{thm}
\label{t-main}
Let $\bal= \{\alpha\ci I\}\ci{I\in\cD}$, $\alpha\ci I\ge 0$, and let $\mu $ and $\nu$ be Radon measures in $\R^N$. Let $1< p <\infty$, and let $1/p + 1/p'=1$. 

The following statements are equivalent
\begin{enumerate}
	\item The following bilinear embedding theorem holds:
\begin{align}
\label{BE-1}
	\sum_{I\in\cD} \left| \int_I f d\mu\right|\cdot \left| \int_I g d\nu\right| \alpha\ci I \le C_1 \|f\|_{L^p(\mu)}  \|g\|_{L^{p'}(\nu)} \qquad \forall f\in L^p(\mu), g\in L^{p'}(\nu). 
\end{align}

\item For all $I_0\in\cD$ 
\begin{align}
\label{BE-test1}
\int_{I_0} \Biggl| \sum_{I\in\cD:\,I\subset I_0} \alpha\ci I \mu(I) \1\ci I  \Biggr|^p d\nu & \le C_2^{p} \mu(I_0) \\
\label{BE-test1'}
\int_{I_0} \Biggl| \sum_{I\in\cD:\,I\subset I_0} \alpha\ci I \nu(I) \1\ci I \Biggr|^{p'} d\mu & \le C_2^{p'} \nu(I_0)
\end{align}

\end{enumerate}
Moreover, $C_2 \le C_1 \le C(p) C_2$
\end{thm}

The condition \cond1 of the theorem equivalent to the fact that the operator $T_\alpha$, $T_\alpha f = \sum_{I\in\cD} \alpha\ci I \left( \int_I f \,d\mu \right) \1\ci I$ is a bounded operator acting from $L^2(\mu)$ to $L^2(\nu)$. 

Condition \cond2 of the theorem is just relaxation of the testing condition $\|T_\alpha \1\ci{I_0}\|\ci{L^p(\nu)} \le C_2 \|\1\ci{I_0}\|\ci{L^p(\mu)}$ and its dual $\|T_\alpha^* \1\ci{I_0}\|\ci{L^{p'}(\mu)} \le C_2 \|\1\ci{I_0}\|\ci{L^{p'}(\nu)}$. Thus the implication \cond1 $\implies$ \cond2 and the estimate $C_2\le C_1$ are trivial. 

The non-trivial part \cond2 $\implies$ \cond1 with the estimate $C_1 \le C(p) C_2$ follows immediately from the proposition below. 

\begin{prop}
\label{p-BE_half}
Let $\cL\subset \cD$ be a collection of dyadic cubes in $\R^N$, and let $f\ge0$ and $g\ge 0$ be functions on $\R^N$ such  that for all $I\in \cL$
\begin{align}
\label{split-cond1}
\mu(I)^{1-p} \left( \int_I f  d\mu \right)^p\ge \nu(I)^{1-p'} \left( \int_I g d\nu  \right)^{p'}.
\end{align}
Let $\alpha\ci I\ge 0$, $I\in\cL$, be such that that for all $I_0\in \cL$ (equivalently for all $I_0\in\cD$)
\begin{align}
\label{test-cond1}
\int_{I_0} \Biggl( \sum_{I\in\cL:\,I\subset I_0} \alpha\ci I \mu(I) \1\ci I  \Biggr)^p d\nu   \le  \mu(I_0) 
\end{align}
Then 
\begin{align*}
\sum_{I\in\cL} \left| \int_I f d\mu\right|\cdot \left| \int_I g d\nu\right| \alpha\ci I \le A \|f\|\ci{L^p(\mu)} \|g\|\ci{L^{p'}(\nu)} + B  \|f\|\ci{L^p(\mu)}^p  , 
\end{align*}
where $A= 2^{1+1/p} p'$, $B= 4(p')^p$. 
\end{prop}

To show that this proposition implies Theorem \ref{t-main} take $f, g$, $\|f\|\ci{L^p(\mu)} = \|g\|\ci{L^{p'}(\nu)} =1$. By Proposition \ref{p-BE_half} the  condition \eqref{BE-test1} implies that the sum over cubes satisfying \eqref{split-cond1} is bounded (by $A+B$). The dual condition \eqref{BE-test1'} implies the estimate of the sum over the rest of the cubes, so Theorem \ref{t-main} is proved for $f, g$, $\|f\|\ci{L^p(\mu)} = \|g\|\ci{L^{p'}(\nu)} =1$. 
The rest follows from the homogeneity. 

\section{Proof of Proposition \ref{p-BE_half} }

\subsection{Stopping moments}
 

Let us 
apply the standard construction of stopping moments (stopping cubes) to  construct the collection $\cG\subset\cL \subset\cD$ of stopping cubes  as follows.%
\footnote{Recall that $\cL\subset \cD$ is the collection of cubes from Proposition \ref{eq-max-norm}. However, the construction works for arbitrary $\cL\subset \cD$.} %
   For a cube $J$ let $\cG^*(I)$ be the collection of maximal cubes $I\in\cL$, $I\subset J$ such that 
\begin{align*}
\mu(I)^{-1} \int_I f d\mu \ge 2 \mu(J)^{-1} \int_J f d\mu
\end{align*}

Let $\cL(J):=\{I\in \cL: I\subset J\}$, and let $G(J) := \cup_{I\in\cG^*(J)} I$. Define also 
\begin{align}
\label{cE(J)}
\cE (J):=\cL(J) \setminus \cup_{K\in \cG^*(J)} \cL(K).
\end{align}

Then it is easy to see that the collection of stopping cubes $\cG^*(I)$ satisfies the following properties:

\begin{enumerate}
	\item For any $I\in \cE(J):= \cL(J) \setminus \cup_{K\in \cG^*(J)} \cL(K)$ we have
	\[
	\mu(I)^{-1} \int_I f\,d\mu < 2  \mu(J)^{-1} \int_J f d\mu  .
	\]
	

\item $\mu(G(J) ) \le \mu(J)/2$.

\end{enumerate}

To construct the collection $\cG$ of stopping cubes, fix some large integer $R$, and consider all maximal $J\in \cL$, $\ell(J)\le 2^R$; that will be the first generation $\cG^*_1$ of stopping cubes. To get the second generation of stopping moments  for each $I\in\cG^*_1$ we  construct the collection $\cG^*(I)$ of stopping moments, and define the second generation $\cG^*_2= \cup_{I\in\cG^*_1}\cG^*(I)$. The next generations are defined inductively,   
\begin{align*}
\cG^*_{n+1}:= \bigcup_{I\in\cG^*_n} \cG^*(I),  
\end{align*}
and we define the collection of stopping cubes $\cG$ by $\cG:= \cup_{n\ge1} \cG^*_n$. 



Propery \cond2 implies that the collection $\cG$ of the stopping cubes satisfies 
the folloowing \emph{Carleson measure condition}
\begin{equation}
\label{CMC}
\sum_{I\in\cG, I\subset J} \mu(I) \le 2 \mu(J) \qquad \forall J\in \cD. 
\end{equation}

We will use the following well-known result. 
 
\begin{lm}[Martingale Carleson Embedding Theorem]
\label{l-dCET}
Let $\mu$ be a measure (on $\R^d$) and let $w\ci I\ge 0$, $I\in\cD$ satisfy the Carleson measure condition
\begin{align}
\label{eq-carl}
\sum_{I\in\cD:\, I\subset J} w\ci I \le C \mu(J) . 
\end{align}
Then for any measurable $f\ge 0$ and for any $p\in(1, \infty)$
\begin{align*}
\sum_{I\in\cD }  \left( \mu(I)^{-1} \int_I f\, d\mu \right)^{p} w\ci I \le (p')^p C \|f\|^p\ci{L^p(\mu)}
\end{align*}
\end{lm}
This lemma (with some constant $C(p)$ instead of $(p')^p$) is well-known. We will explain the constant $(p')^p$  later. 

\subsection{Splitting the estimate}

%

Since the collection of cubes $I\in\cL$ such that $\ell(I)\le 2^K$ can be represented as the union $\cup_{J\in\cG} \cE(J) $  we can write
\begin{align*}
 \sum_{I\in\cL: \ell(I)\le 2^R} \left( \int_I f d\mu\right)\cdot \left( \int_I g d\nu\right) \alpha\ci I     
 \le \sum_{J\in  \cG} \sum_{I\in\cE(J)}    \left( \int_I f d\mu\right)\cdot \left( \int_I g d\nu\right) \alpha\ci I 
\end{align*}
We can represent the inner sum as an integral
\begin{align}
\label{int-sum_EJ}
\sum_{I\in\cE(J)}    \left( \int_I f d\mu\right)\cdot \left( \int_I g d\nu\right) \alpha\ci I  
= \int F\ci{\cE(J)} g \, d\nu , 
\end{align}
where 
\[
F\ci{\cE(J)}:=  \sum_{I\in\cE(J) }\alpha\ci I  \left( \int_I f d\mu\right)  \1\ci I
\]

The above property \cond1 of $\cE(J)$ 
imply that for $I\in\cE(J)$
\[
  \int_I f d\mu < 2 \mu(I) \left(\mu(J)^{-1} \int_J f d \mu  \right) .
\]
Then the condition \eqref{test-cond1} implies that 
\begin{align}
\label{Lp-norm_F_EJ}
\|F\ci{\cE(J)}\|\ci{L^p(\nu)}^p \le 2^p \left(\mu(J)^{-1} \int_J f d \mu  \right)^p \mu(J)  .
\end{align}

We now split the integral in \eqref{int-sum_EJ}, 
\[
\int F\ci{\cE(J)} g \, d\nu  = \int_J F\ci{\cE(J)} g \, d\nu = \int_{J\setminus G(J)} F\ci{\cE(J)} g \, d\nu 
+ \int_{G(J)} F\ci{\cE(J)} g \, d\nu  = A(J) + B(J). 
\]
The main reason for this splitting is that the sets $J\setminus \cG(J)$, $J\in \cG$ are disjoint, so the sum of $A(J)$ is easy to estimate. 

\subsection{The easy estimate}
The sum of $A(J)$ is easy to estimate.  
Namely, using \eqref{Lp-norm_F_EJ} we can write
\begin{align*}
\sum_{J\in\cG} A(J)  & \le  \sum_{J\in\cG} \| F\ci{\cE(J)} \|\ci{L^p(\nu)} \| g \1\ci{J\setminus G(J)} \|\ci{L^{p'}(\nu)} 
\\ & \le
\left(  \sum_{J\in\cG} \| F\ci{\cE(J)} \|\ci{L^p(\nu)}^p \right)^{1/p} 
\left(  \sum_{J\in\cG}   \| g \1\ci{J\setminus G(J)} \|\ci{L^{p'}(\nu)}^{p'}  \right)^{1/p'}  && \text{H\"{o}lder inequality}
\\  &  \le 
\left(  \sum_{J\in\cG} \| F\ci{\cE(J)} \|\ci{L^p(\nu)}^p \right)^{1/p}  \| g  \|\ci{L^{p'}(\nu)} && 
J\setminus G(J)\ \text{are disjoint}
\\  & \le
2\left(  \sum_{J\in\cG}  \left(\mu(J)^{-1} \int_J f d \mu  \right)^p \mu(J)\right)^{1/p}  \| g  \|\ci{L^{p'}(\nu)} && \text{by \eqref{Lp-norm_F_EJ} } 
\end{align*}
Applying  Lemma \ref{l-dCET} with 
\[
w\ci I = \left\{
\begin{array}{ll} \mu(I) , \qquad & I \in \cG \\ 0 & I\notin \cG .\end{array} \right. 
\]
we get using the Carleson measure property \eqref{CMC}  that 
\begin{equation}
\label{C-embed-1}
\sum_{J\in\cG}  \left(\mu(J)^{-1} \int_J f d \mu  \right)^p \mu(J) \le 2 (p')^p \| f\|\ci{L^p(\mu)}^p, 
\end{equation}
so 
%
\[ 
\sum_{J\in \cD} A(J) \le  2^{1+1/p} p'
\|f\|\ci{L^p(\mu)} \| g\|\ci{L^{p'}(\nu)} .
\]

\subsection{``Replacing the averages'' and the ``hard'' estimate}  Let us now estimate $\sum_{J\in \cL} B(J)$. This is the part where we use the splitting condition \eqref{split-cond1}. 

Recall that 
\begin{align}
\label{eq-BJ}
B(J) = \int_{G(J)} F\ci{\cE(J)} g \, d\nu 
\end{align}
and that $G(J) = \bigcup_{I\in\cG^*(J)} I$. Since $F\ci{\cE(J)}$ is constant  on the intervals $I\in \cG^*(J)$, one can replace $g$ in \eqref{eq-BJ} by the function 
\[
\wt g\ci J := \sum_{I\in\cG^*(J)} \left( \nu(I)^{-1} \int_I g\, d\nu \right) \1\ci I .
\]
Then one can estimate
\begin{align*}
B(J)  & = \int_J F\ci{\cE(J)} \wt g\ci J    \le \| F\ci{\cE(J)} \|\ci{L^p(\nu)} \| \wt g\ci J \|\ci{L^{p'}(\nu)}
&&
\\
& \le   2 \left(\mu(J)^{-1} \int_J f d \mu  \right) \mu(J)^{1/p} \| \wt g\ci J \|\ci{L^{p'}(\nu)} && \text{by \eqref{Lp-norm_F_EJ} }
\\
& = 
2 \left(\mu(J)^{-1} \int_J f d \mu  \right) \mu(J)^{1/p}
\left(  
\sum_{I\in \cG^*(J)} \left( \nu(I)^{-1} \int_I g\, d\nu \right)^{p'} \nu(I) 
\right)^{1/p'}
\\
& 
\le 2 \left(\mu(J)^{-1} \int_J f d \mu  \right) \mu(J)^{1/p}
\left(  
\sum_{I\in \cG^*(J)} \left( \mu(I)^{-1} \int_I f\, d\mu \right)^{p} \mu(I) 
\right)^{1/p'}
&& \text{by \eqref{split-cond1} }
\end{align*}
Therefore, summing over all generations of stopping cubes and using H\"{o}lder inequality we get
\begin{align*}
\sum_{J\in \cG} B(J) 
& \le 
2\left(\sum_{J\in\cG}  \left(\mu(J)^{-1} \int_J f d \mu  \right)^p \mu(J) \right)^{1/p} 
\left( \sum_{J\in\cG} \sum_{I\in\cG^*(J)} \left( \mu(I)^{-1} \int_I f\, d\mu \right)^{p} \mu(I) \right)^{1/p'} 
\\
& \le 2 \sum_{J\in\cG}  \left(\mu(J)^{-1} \int_J f d \mu  \right)^p \mu(J) ;
\end{align*}
the last inequality holds because the sum in the second term is the sum over all $I\in \cG\setminus \cG_1$, (where, recall, $\cG_1$ is the first generation of stopping cubes) so it is dominated by the sum in the first term. 

But the final sum was already estimated in \eqref{C-embed-1}! So
\[
\sum_{J\in\cG} B(J)\le 4 (p')^p \|f\|\ci{L^p(\mu)}^p 
\]

%

\subsection{Concluding the proof} Gathering all the estimates together we get
\begin{align*}
 \sum_{I\in\cL: \ell(I)\le 2^R} \left( \int_I f d\mu\right)\cdot \left( \int_I g d\nu\right) \alpha\ci I  
 &\le \sum_{J\in\cG} A(J)
 +  \sum_{J\in\cG} B(J) 
 \\
 &\le 2^{1+1/p} p' \|f\|\ci{L^p(\mu)} \|g\|\ci{L^{p'}(\nu)} + 4 (p')^p \|f\|\ci{L^p(\mu)}^p, 
\end{align*}
and the right side does not depend on $R$. Letting $R\to\infty$ we get the conclusion of the proposition. \hfill \qed

\section{Maximal function and dyadic Carleson Embedding theorem}
This section contains well-known facts and is presented only to save a reader a trip to a library. We give here a quick explanation of why Lemma \ref{l-dCET} holds with the constant $(p')^p$. 

One of the standard way of proving the Carleson Embedding type of result is a the comparison with the maximal function. Recall that given a Radon measure $\mu$ in $\R^d$ the dyadic maximal function $M_\mu = M_\mu^{\text{d}}$ is defined by 
\begin{align*}
M_\mu f (x) = \sup_{I\in\cD:\, I \ni x} \mu(I)^{-1} \left| \int_I f(x) \, dx \right|
\end{align*}
The maximal function operator $M_\mu$ is a particular case of a martingale maximal function (when one restricts everything to a finite cube), so Theorem 14.1 from \cite{Wil-ProbMart-1991} which states that martingale maximal function is bounded in $L^p$, $p\in(1, \infty)$ with the norm at most $p'$ implies that 
\begin{align}
\label{eq-max-norm}
\|M_\mu f \|\ci{L^p(\mu)} \le p' \|f\|\ci{L^p(\mu)}. 
\end{align}

The Carleson Embedding Theorem can be obtained from this result by the standard level sets comparison. Namely, for $\lambda>0$ let $\cE_\lambda$ be the collection of cubes $Q$ such that
\begin{align*}
\mu(Q)^{-1} \left| \int_Q f(x) \, dx \right| >\lambda , 
\end{align*}
and let $E_\lambda := \bigcup_{Q\in \cE_\lambda} \cE_\lambda$. Then clearly $M_\mu f (x) > \lambda$ on $E_\lambda$, i.e. $E_\lambda \subset \{ x\in \R^d\,:\, M_\mu f(x)>\la \}$. 

On the other hand, since the set $E_\la$ can be represented as a disjoint union  of maximal cubes in $\cE_\la$, condition \eqref{eq-carl} of Lemma \ref{l-dCET} implies that 
\[
\sum_{Q\in \cE_\la} w\ci Q \le C \mu( E_\la ). 
\]
But since $E_\la$ is contained in the sublevel set of $M_\mu f$
\[
\mu( E_\la ) \le \mu \left(  \{ x\in \R^d\,:\, M_\mu f(x)>\la \}  \right), 
\]
so 
\[
\sum_{Q\in \cE_\la} w\ci Q \le C \mu \left(  \{ x\in \R^d\,:\, M_\mu f(x)>\la \}  \right).
\]
Therefore
\[
\sum_{Q\in\cD} \left( \mu(Q)^{-1} \int_Q f d\mu \right)^p w\ci Q \le 
C \int_{\R^d} ( M_\mu f)^p d\mu \le C (p')^p \|f\|\ci{L^p(\mu)}^p;
\]
the last inequality here follows from \eqref{eq-max-norm}. \hfill\qed

\

\section{A concluding remark}

As a reader could see, the above construction used none of the specific properties of the dyadic lattice $\cD$. In fact, all the proofs work in a more general martingale situation. 

Namely, one can consider a set $\cX$ with two $\sigma$-finite measures $\mu$ and $\nu$ (both defined on the same $\sigma$-algebra $\fA$) and a collection 
(lattice)  
$\cD=\cup_{k\in\Z} \cD_k$ of $\fA$-measurable sets, such that for each $k$ the collection $\cD_k$ is a countable partition of $\cX$ and $\cD_{k+1}$ is a refinement of $\cD_k$. 

All the proofs work in this setting, one literally does not have to change anything.  




\end{document}